\documentclass[graybox]{svmult}

\usepackage{mathptmx}       % selects Times Roman as basic font
\usepackage{helvet}         % selects Helvetica as sans-serif font
\usepackage{courier}        % selects Courier as typewriter font
\usepackage{type1cm}        % activate if the above 3 fonts are
\usepackage{makeidx}         % allows index generation
\usepackage{graphicx}        % standard LaTeX graphics tool
\usepackage{multicol}        % used for the two-column index
\usepackage[bottom]{footmisc}% places footnotes at page bottom

\usepackage{etex} 
\reserveinserts{30}
\usepackage{amsmath, amssymb}
\usepackage{graphicx}
\usepackage{tikz}
\usepackage{xcolor,calc}
\usetikzlibrary{matrix,arrows,decorations.pathmorphing}
\usepackage{amsmath,amscd}%
\usepackage{amsfonts}%
\usepackage{amssymb}%
\usepackage{graphicx}
\usepackage{txfonts}
\usepackage{color}
\usepackage{subfigure}
\usepackage{pgfplots}
\usepackage[sort&compress]{natbib}
\makeindex             % used for the subject index
\newcommand{\kdifform}[2]{#1^{(#2)}} % differential k-forms with rank
\newcommand{\incidencederivative}[2]{\mathsf{E}^{(#1,#2)}} % Incidence matrix on k-cochains, implements the  exterior derivative operator, or dual of the boundary operator
\newcommand{\innerspace}[3]{\left(#1,#2\right)_{#3}} % inner product integrated

\newcommand{\inner}[2]{\left( #1, #2\right)} % inner product

\newcommand{\projection}{\pi_{h}} % projection operator
\newcommand{\ederiv}{\mathrm{d}} % exterior derivative continuous
\newcommand{\spacemap}{\rightarrow} % arrow that maps from one space to another
 
\newcommand{\figref}[1]{Figure~\ref{#1}} % reference a figure
\newcommand{\realNumber}{\mathbb{R}}
\newcommand{\differentialFormSpace}[2]{\Lambda^{#1}\left(#2\right)}

\newcommand{\differentialFormSpaceDiscrete}[2]{\Lambda^{#1}_{h}\left(#2\right)}

\begin{document}

\title*{Mimetic Spectral Element Advection}
\author{Artur Palha, Pedro Pinto Rebelo and Marc Gerritsma}
\institute{Artur Palha, Pedro Pinto Rebelo and Marc Gerritsma \at Delft University of Technology, Faculty of Aerospace Engineering, Aerodynamics Group, 2600GB - Delft,The Netherlands, \email{\{A.Palha, P.J.PintoRebelo, M.I.Gerritsma\}@tudelft.nl}}
\maketitle

\abstract{We present a discretization of the linear advection of differential forms on bounded domains. The framework established in \cite{kreeft2011mimetic} is extended to incorporate the Lie derivative, $\mathcal L$, by means of Cartan's homotopy formula. The method is based on a physics-compatible discretization with spectral accuracy . It will be shown that the derived scheme has spectral convergence with local mass conservation. Artificial dispersion depends on the order of time integration.}

\abstract*{We present a discretization of the linear advection of differential forms on bounded domains. The framework established in \cite{kreeft2011mimetic} is extended to incorporate the Lie derivative, $\mathcal L$, by means of Cartan's homotopy formula. The method is based on a physics-compatible discretization with spectral accuracy . It will be shown that the derived scheme has spectral convergence with local mass conservation. Artificial dispersion depends on the order of time integration.}

\section{INTRODUCTION}
\label{Section::Introduction}

Consider the classical advection problem for a scalar function in conservation form,
\begin{equation}
 \frac{\partial \rho}{\partial t} + \nabla\cdot\left(\vec{v}\rho\right) = 0,
 \label{eq::scalar_advection_diffusion}
\end{equation}
where $\vec{v}$ is a prescribed  uniformly Lipschitz continuous vector field and $\rho$ the advected scalar function. The method presented in this work is based on the approximation of the differential operators with the focus on the spatial discretization and a time-stepping scheme that distinguishes between quantities evaluated at time instants and quantities evaluated over time intervals. 
%Although the method presented here will be applied only to the advection of scalar quantities, \eqref{eq::scalar_advection_diffusion}, it can be extended to the non-scalar case. 

The mimetic framework, presented in \cite{kreeft2011mimetic}, showed that using the differential geometric approach for the representation of physical laws clarifies the underlying structures. One clearly identifies to what kind of geometrical object a certain physical quantity is associated and this determines how its discretization must be done (e.g.: evaluation at points, integration over lines, surfaces or volumes). Additionally, a well defined, metric free, representation of differential operators is obtained, together with their metric dependent Hilbert adjoints. For these reasons, the authors followed this approach for the advection equation. It is known, see \cite[pp. 317]{abraham_diff_geom}, that \eqref{eq::scalar_advection_diffusion} is a particular case of the generalized advection equation which can be written in terms of differential geometry as,
\begin{equation}
 \frac{\partial\kdifform{\alpha}{k}}{\partial t} + \mathcal{L}_{\vec{v}}\, \kdifform{\alpha}{k} = 0 \label{eq::testCases_advection_diffusion}.
\end{equation} 
The advection operator, $\mathcal{L}_{\vec{v}}$, is the Lie derivative for the prescribed velocity field $\vec{v}$ and the advected quantity is given by the $k$-differential form $\kdifform{\alpha}{k}$. Depending on the index $k$ the quantity $\alpha^{(k)}$ can represent scalar,  vector and higher dimensional quantities.

\section{DIFFERENTIAL GEOMETRY}
\label{Section::DifferentialGeometry}

In this section a brief introduction to differential geometry is given. For a more detailed introduction the reader is directed to \cite{abraham_diff_geom}.
Given an $n$-dimensional smooth orientable manifold $\Omega$ it is possible to define in each point a tangent vector space $E$ of dimension $n$. The space of smooth vector fields on a manifold is the space, $\Gamma$, of smooth assignments of elements of $E$ to each point of the manifold. We denote by $\Lambda^{k}$,  $k$ an integer $0 \leq k \leq n$, the space of differentiable $k$-forms, i.e. all smooth $k$-linear, antisymmetric maps $\kdifform{\omega}{k}: E \times \cdots \times E \rightarrow \realNumber$, at every point of the manifold. We recall the wedge product $\wedge : \Lambda^{k} \times \Lambda^{l} \rightarrow \Lambda^{k+l}$ for $k+l \leq n$ with the property that $\kdifform{\alpha}{k} \wedge \kdifform{\beta}{l} = (-1)^{kl} \kdifform{\beta}{l} \wedge \kdifform{\alpha}{k}$. The inner product $\inner{\cdot}{\cdot}$ on $E$ induces at each point of the manifold an inner product $\inner{\cdot}{\cdot}$ on $\Lambda^{1}$. In turn, this can be extended to a local inner product on $\Lambda^{k}$, \cite[pp. 149]{morita::geometry_diff_forms}. The local inner product gives rise to a unique metric operator, Hodge-$\star$, $\star: \Lambda^{k} \rightarrow \Lambda^{n-k}$, defined by $\kdifform{\alpha}{k} \wedge \star \kdifform{\beta}{k} = \innerspace{\kdifform{\alpha}{k}}{\kdifform{\beta}{k}}{}\omega^{(n)}$, where $\omega^{(n)}=\star 1$ is the standard volume form. By integration, one can define an inner product on $\Omega$ as $\inner{\cdot}{\cdot}_{L^{2}}:=\int_{\Omega}\inner{\cdot}{\cdot} \omega^{(n)}$. The exterior derivative $\ederiv : \Lambda^{k} \rightarrow \Lambda^{k+1}$ satisfies the following rule, $\ederiv\left(\kdifform{\alpha}{k} \wedge \kdifform{\beta}{l} \right) = \ederiv \kdifform{\alpha}{k} \wedge \kdifform{\beta}{l} +(-1)^k \kdifform{\alpha}{k} \wedge \ederiv \kdifform{\beta}{l}$ and by definition $\ederiv \alpha^{(n)} = 0$. The flat operator, $\flat$, is a mapping $\flat:\Gamma\mapsto\Lambda^{1}$. 

The Lie derivative along a tangent vector field, $\vec{v}$, is denoted by $\mathcal{L}_{\vec{v}}$ and represents the advection operator in differential geometry. It is a mapping $\mathcal{L}_{\vec{v}}:\Lambda^{k}\mapsto\Lambda^{k}$. From Cartan's homotopy formula the Lie derivative can be written as
\[
\mathcal{L}_{\vec{v}}\, \kdifform{\alpha}{k} ():= \ederiv \iota_{\vec{v}}\kdifform{\alpha}{k} + \iota_{\vec{v}} \ederiv \kdifform{\alpha}{k}\;,
\]

 where the interior product of a tangent vector field, $\vec{v}$, with a $k$-form, $\kdifform{\alpha}{k}$, is a mapping $\iota_{\vec{v}}\kdifform{\alpha}{k}: \Lambda^{k} \rightarrow \Lambda^{k-1}$ given by:
\[
	\iota_{\vec{v}}\kdifform{\alpha}{k}(\vec{X}_{2},\cdots,\vec{X}_{k}) := \kdifform{\alpha}{k}\left(\vec{v},\vec{X}_{2}, \cdots,\vec{X}_{k} \right),\quad \forall \vec{X}_{i}\in\Gamma\quad\mathrm{and}\quad\iota_{\vec{v}}\kdifform{\alpha}{0}=0,\quad\forall\vec{v}\in\Gamma\;.
\]
\begin{svgraybox}
The interior product is the adjoint of the wedge product,  made explicit by:
\begin{equation}
\inner{\iota_{\vec{v}}\,\kdifform{\alpha}{k}}{\kdifform{\beta}{k-1}}_{L^{2}\Lambda^{k-1}} = \inner{\kdifform{\alpha}{k}}{\vec{v}^{\flat}\wedge\kdifform{\beta}{k-1}}_{L^{2}\Lambda^{k}},\quad\forall\kdifform{\beta}{k-1}\in\Lambda^{k-1}\label{eq::adjoint_interior_product}
\end{equation}
where $\vec{v}^{\,\flat}=\kdifform{\nu}{1}\in\Lambda^{1}$ and $\kdifform{\alpha}{k}_{h}\in \Lambda^{k}$.
\end{svgraybox}

The relevance of this adjoint relation between the interior product and the wedge product lies in the fact that it shows how a physical quantity represented by an interior product with a vector field can be represented by its dual differential 1-form.

For a volume form $\kdifform{\rho}{n}$ the Lie derivative is simply $\mathcal{L}_{\vec{v}}\, \kdifform{\rho}{n} =  \ederiv \iota_{\vec{v}}\kdifform{\rho}{n}$ and for a $0$-form, $\mathcal{L}_{\vec{v}}\, \kdifform{\alpha}{0} =  \iota_{\vec{v}} \ederiv \kdifform{\alpha}{0}$.

\section{MIMETIC DISCRETIZATION}
\label{Section::Discretization}

In this section a brief introduction to the discretization of physical quantities and to the discretization of the exterior derivative is presented. For a more detailed presentation the reader is directed to \cite{kreeft2011mimetic,palhaLaplaceDualGrid,gerritsma::geometricBasis}

Consider a three dimensional domain $\Omega$ and an associated grid consisting of a collection of points, $\tau_{(0),i}$, line segments connecting the points, $\tau_{(1),i}$, surfaces bounded by these line segments, $\tau_{(2),i}$, and volumes bounded by these surfaces, $\tau_{(3),i}$.

\begin{svgraybox}
Let $\Lambda^{k}$ be the space of smooth differentiable $k$-forms. Additionally, let the finite dimensional space of differentiable forms be defined as $\Lambda^{k}_{h}=\mathrm{span}(\{\kdifform{\epsilon_{i}}{k}\}), \,\, i=1,\cdots,\dim(\Lambda^{k}_{h})$, where $\kdifform{\epsilon_{i}}{k}\in \Lambda^{k}$ are basis $k$-forms. Under these conditions it is possible, see \cite{kreeft2011mimetic,palhaLaplaceDualGrid}, to define a projection operator $\pi_{h}$ which projects elements of $\Lambda^{k}$ onto elements of $\Lambda^{k}_{h}$ which satisfies:
	\begin{equation}
		\pi_{h}\ederiv = \ederiv \pi_{h}\;.  \label{eq:commuting_projection}
	\end{equation}
	It is possible to write:
	\[
		\pi_{h} \kdifform{\alpha}{k} = \kdifform{\alpha_{h}}{k}=\sum_{i}\alpha_{i}\kdifform{\epsilon_{i}}{k}\;,
	\]
	where
	\[
		\alpha_{i} = \int_{\tau_{(k),i}}\kdifform{\alpha}{k} \quad\mathrm{and}\quad  \int_{\tau_{(k),i}}\kdifform{\epsilon_{j}}{k} = \delta_{ij},\quad k=0,1,\cdots,n\;.
	\]
\end{svgraybox}
A set of basis functions yielding a projection operator $\pi_{h}$ that satisfies \eqref{eq:commuting_projection} can be constructed using piecewise polynomial expansions on the quadrilateral elements using tensor products. Thus, it suffices to derive the basis forms in one dimension on a reference interval and generalize them in $n$ dimensions. 

In one dimension take a 0-form, $\alpha^{(0)} \in \differentialFormSpace{0}{Q_{ref}}$, where $Q_{ref} := \left[ -1,1 \right]$. Define on $Q_{ref}$ a cell complex $D$ of order $p$ consisting of $(p+1)$ nodes $\tau_{(0),i}=\xi_{i}$ with $i=0,\cdots,p$, where $-1 \leq \xi_{0} < \cdots < \xi_{i}< \cdots \xi_{p}\leq 1$ are the Gauss-Lobatto quadrature nodes, and $p$ edges, $\tau_{(1),i} = \left[ \xi_{i-1},\xi_{i} \right]$ with $i=1,\cdots,p$. The projection operator $\pi_{h}$ reads:
\begin{align}
\projection \alpha^{(0)} \left( \xi \right) = \sum_{i=0}^{p} \alpha_{i} \epsilon^{(0)}_{i}(\xi)\;,
\end{align}
where $\epsilon^{(0)}_{i}(\xi) = l_{i} \left( \xi \right)$ are the $p^{th}$ order \emph{Lagrange polynomials} and $\alpha_{i}=\alpha^{0}(\xi_{i})$.
 Similarly in one dimension for the projection of 1-forms Gerritsma \cite{gerritsma::edge_basis} and Robidoux \cite{robidoux-polynomial} derived 1-form polynomials called \emph{edge polynomials}, $\epsilon^{(1)}_{i} \in \differentialFormSpaceDiscrete{1}{Q_{ref}}$,
\begin{align}
\epsilon^{(1)}_{i} (\xi) = e_{i} \left( \xi \right) \ederiv\xi, \quad \text{with} \quad e_{i}(\xi) = - \sum_{k=0}^{i-1} \frac{d l_{k}}{d \xi}\;.
\end{align}
Note that in this way we have:
\begin{align}
\int_{\xi_{j-1}}^{\xi_{j}}\epsilon^{(1)}_{i} = \int_{\xi_{j-1}}^{\xi_{j}} e_{i} \left( \xi \right)\ederiv\xi = \delta_{ij}\;.
\end{align}

Moreover, the exterior derivative of the basis 0-forms is given by:
	\begin{equation}
		\ederiv\kdifform{\epsilon_{i}}{0} = \frac{\ederiv l_{i}}{\ederiv\xi}\,\ederiv\xi = -\sum_{k=0}^{i-1}\frac{\ederiv l_{k}}{\ederiv \xi}\,\ederiv\xi - \left(-\sum_{k=0}^{i}\frac{\ederiv l_{k}}{\ederiv \xi}\,\ederiv\xi\right) = \kdifform{\epsilon_{i}}{1} -\kdifform{\epsilon_{i+1}}{1},\quad i=1,\cdots,p-1\;.
	\end{equation}
	In this way, the exterior derivative of a discrete 0-form can be written as:
	\begin{equation}
		\ederiv\kdifform{\alpha_{h}}{0} = \ederiv\sum_{i=0}^{p}\alpha_{i}\kdifform{\epsilon_{i}}{0} = \sum_{i=0,j=1}^{p}\incidencederivative{1}{0}_{ij}\alpha_{j}\kdifform{\epsilon_{i}}{1}\;, \label{incidence_d}
	\end{equation}
	where, $\incidencederivative{1}{0}_{ij}$ is the incidence matrix containing only the values, 0, 1 and -1, see \cite{kreeft2011mimetic,palhaLaplaceDualGrid,gerritsma::geometricBasis} for more details. This idea can be extended to higher dimensions, giving rise to $k$-incidence matrices, $\incidencederivative{k+1}{k}_{ij}$, which represent the discrete exterior derivative on discrete $k$-forms, see \cite{palhaLaplaceDualGrid,gerritsma::geometricBasis}.

%The projection operator defined in this way in one dimension satisfies \eqref{eq:commutation_relations_projection_derivative}, see \cite{kreeft2011mimetic}. For higher dimensional spaces tensor products are used. For example, for the two dimensional case we have:
%\begin{equation}
%	\kdifform{\epsilon_{ij;p}}{0} = e^{p}_{i}(\xi)l^{p}_{j}(\eta), \quad 
%	\begin{cases}
%		\kdifform{\epsilon_{ij,\xi;p}}{1} = e^{p}_{i}(\xi)l^{p}_{i}(\eta)\ederiv\xi\\
%		\\
%		\kdifform{\epsilon_{ij,\eta;p}}{1} = l^{p}_{i}(\xi)e^{p}_{i}(\eta)\ederiv\eta
%	\end{cases}
%	,\quad\ \kdifform{\epsilon_{ij;p}}{2} = e^{p}_{i}(\xi)e^{p}_{j}(\eta)\ederiv\xi\ederiv\eta \label{eq::basis_forms}
%\end{equation}
%

\section{Mimetic spectral advection: an application to 1D advection}
\label{Section::Numerical}

In this section we want to illustrate how to discretize the advection equation. Take the Lie advection of a 1-form, 
\begin{equation}
\frac{\partial\kdifform{\rho}{1}}{\partial t} + \ederiv\iota_{\vec{v}}\kdifform{\rho}{1} = 0  \Leftrightarrow\left\{
\begin{array}{l}
\frac{\partial\kdifform{\rho}{1}}{\partial t} = -\ederiv\kdifform{\varsigma}{0}\\
\\
\iota_{\vec{v}}\,\kdifform{\rho}{1} = \kdifform{\varsigma}{0}
\end{array}\;.
\right.  \label{eq::system_advection}
\end{equation}
Here $\kdifform{\rho}{1}$ is the advected quantity, say mass density, and $\kdifform{\varsigma}{0}$ represents the instantaneous fluxes of the advected quantity under the vector field $\vec{v}$, which are discretized in space as,
\begin{equation}
\kdifform{\varsigma_{h}}{0} = \sum_{i=0}^{p}\varsigma_{i}(t)\kdifform{\epsilon_{i}}{0}\quad\mathrm{and}\quad\kdifform{\rho_{h}}{1} = \sum_{i=1}^{p}\rho_{i}(t)\kdifform{\epsilon_{i}}{1}\;. \label{eq::discretization_variables}
\end{equation}
%The system of equations \eqref{eq::system_advection} comprises one equation that defines the evolution in time of the advected quantity, the top equation, and one equation that defines the instantaneous fluxes, the bottom equation. For this reason, and 
For the sake of clarity in the method presentation we first introduce the time treatment, then the interior product discretization and finally their combination for a numerical solution of the advection problem.

\subsection{Time integration}
The time integrator used for solving the time evolution part of the advection equation is the canonical mimetic one, an arbitrary order symplectic operator derived in \cite{palha::phd}, which is connected to canonical Gauss collocation integrators. Take an ordinary differential equation of the unknown function $y(t)$:
\begin{equation}
\frac{\ederiv y}{\ederiv t} = h(y,t), \quad t\in I\subset\mathbb{R}\;.\label{eq:time_ode}
\end{equation}
Discretizing $y(t)$ as $y_{h} = \sum_{k=0}^{p}y^{k}l_{k}(t)$, one gets:
\[
	\frac{\ederiv y_{h}}{\ederiv t} = \sum_{k=0}^{p} y^{k}\frac{\ederiv l_{k}(t)}{\ederiv t} = \sum_{k=1}^{p}(y^{k} - y^{k-1})\,e_{k}(t)\;,
\]
where the superscript $k$ denotes the time level.

\begin{figure}[ht]
\centering
\includegraphics[width=0.4\textwidth]{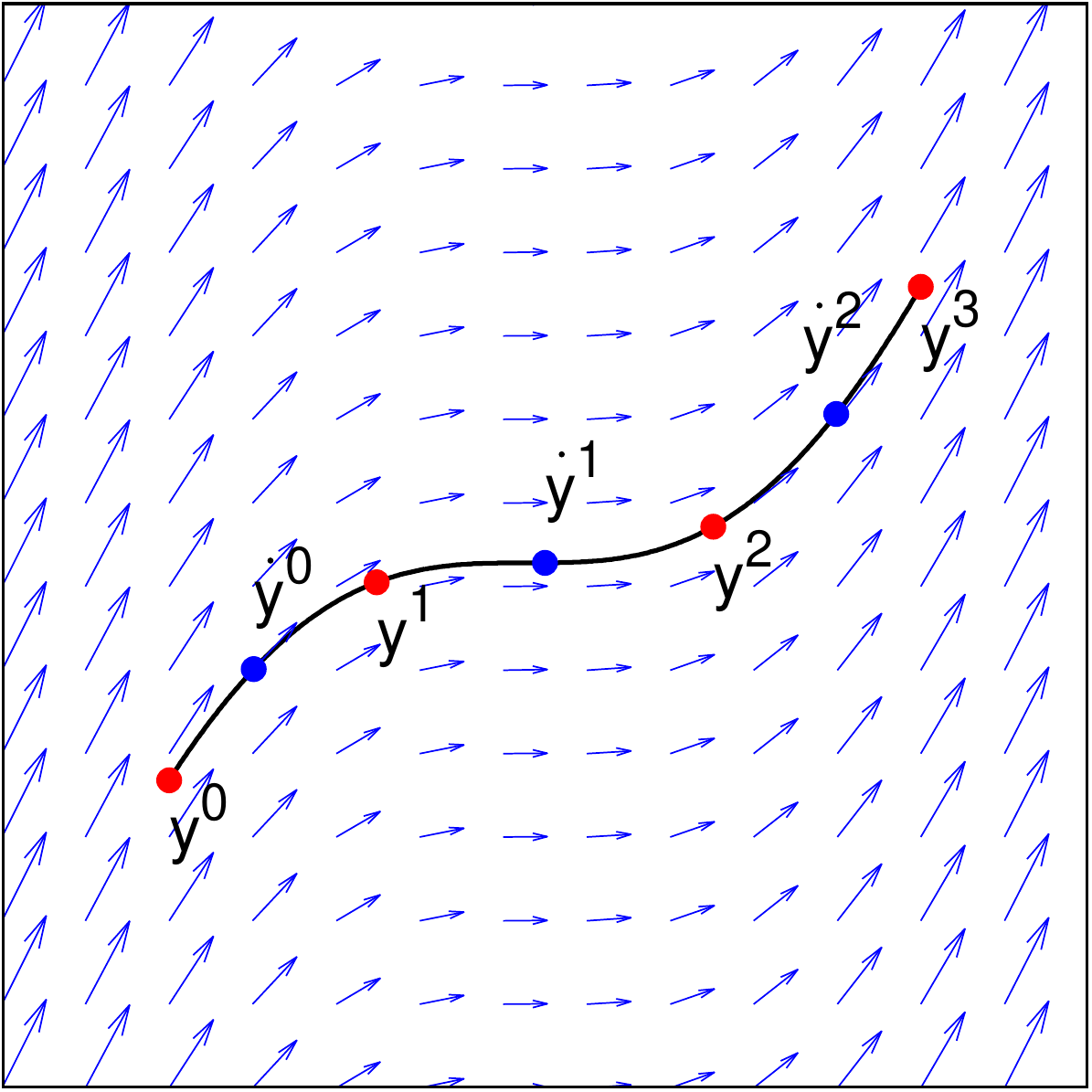}
\caption{Geometric interpretation of the solution of \eqref{eq:time_ode} as given by  \eqref{eq::TestCases_primal_dual_integrator}: $(t,y^{(0)}(t))$. In red the Gauss-Lobatto nodes where the trajectory is discretized. In blue, the Gauss nodes where its derivative is discretized. The flow field, represented by arrows, is tangent to the curve at the Gauss nodes. That is, the derivative of the approximate trajectory is exactly equal to the flow field at the Gauss nodes.}
\label{fig:mimetic_integrator}
\end{figure}

The approximated solution, $y_{h}(t)$, is a polynomial of order $p$ determined by means of $(p+1)$ degrees of freedom such as its values at the Gauss-Lobatto nodes, red dots in \figref{fig:mimetic_integrator}. On the other hand, $\frac{\ederiv y_{h}}{\ederiv t}$ is a polynomial of order $(p-1)$ defined by only $p$ degrees of freedom. One can set these degrees of freedom to be the values of the derivative in one point inside each of the $p$ intervals $[t^{k}, t^{k+1}]$. A choice that results in a symplectic integrator of order $2p$ is to select these points as the Gauss nodes of order $(p-1)$, the blue nodes of \figref{fig:mimetic_integrator}. Notice that  along the trajectory these nodes will not show the usual Gauss-Lobatto and Gauss distribution patterns, since in general the velocity field is not constant. In this way the discrete integrator becomes:
\begin{svgraybox}
\begin{equation}
\sum_{k=1}^{p}(y^{k} - y^{k-1})\,e_{k}(\tilde{t}^{q}) = h\left(\sum_{k=0}^{p}y^{k} l^{k}(\tilde{t}^{q}),\tilde{t}^{q}\right),\quad q=1,2,\cdots,p\;,
\label{eq::TestCases_primal_dual_integrator}
\end{equation}
\end{svgraybox}
with $\tilde{t}^{j}$ the $p$ nodes of a Gauss quadrature formula. The fact that the instants in time, $t^{k}$, where the $y^{k}_{i}$ are defined alternate with the instants in time, $\tilde{t}^{j}$, where the $h_{i}$ are evaluated (see \figref{fig:mimetic_integrator}), corresponds to a staggering in time. This staggering also appears in leap-frog methods and in the implicit midpoint rule, for instance.
 
%\begin{equation}
%\begin{cases}
%\dfrac{\ederiv\sum_{i}\rho_{i}(t)\kdifform{\epsilon_{i}}{1}}{\ederiv t} = -\ederiv\sum_{i}\varsigma_{i}(t)\kdifform{\epsilon_{i}}{0}\\
%\\
%\sum_{i}\rho_{i}(t)\inner{\kdifform{\epsilon_{i}}{1}}{\star\nu^{(0)}\wedge\kdifform{\epsilon_{j}}{1}}_{L^{2}\Lambda^{1}(\Omega_{m})} = \sum_{i}\varsigma_{i}(t)\inner{\kdifform{\epsilon_{i}}{0}}{\kdifform{\epsilon_{j}}{0}}_{L^{2}\Lambda^{0}(\Omega_{m})}\;,\quad\forall\kdifform{\epsilon_{j}}{0}\in H^{1}\Lambda^{0}(\Omega_{m})
%\end{cases}
%\label{eq::system_advection_discretization_02}
%\end{equation}
%where, for simplification, it was assumed that the domain $\Omega$ is Euclidean. Recalling that in the discrete space of differentiable forms the exterior is related to the incidence matrix, \eqref{incidence_d}, one can substitute the action of the exterior derivative on the basis 0-forms by a linear combination of the basis 1-forms. Therefore one obtains for the first equation of \eqref{eq::system_advection_discretization_02}:
%

The first equation in \eqref{eq::system_advection} using the discretization \eqref{eq::discretization_variables} and \eqref{incidence_d} can be written as:
\begin{equation}
\dfrac{\sum_{i}\ederiv\rho_{i}(t)\kdifform{\epsilon_{i}}{1}}{\ederiv t} = -\sum_{il}\incidencederivative{1}{0}_{il}\varsigma_{l}(t)\kdifform{\epsilon_{i}}{1}\quad\Rightarrow\quad\dfrac{\ederiv\rho_{i}(t)}{\ederiv t} = -\sum_{l}\incidencederivative{1}{0}_{il}\varsigma_{l}(t)\;.
\label{eq::system_advection_discretization_03}
\end{equation}
%For this equality to hold, one needs to satisfy:
%\begin{equation}
%\dfrac{\ederiv\rho_{i}(t)}{\ederiv t} = -\sum_{l}\incidencederivative{1}{0}_{il}\varsigma_{l}(t)\;.
%\label{eq::system_advection_discretization_04}
%\end{equation}
This equation has a similar form as \eqref{eq:time_ode}, but now as a system of equations, therefore one can apply the mimetic integrator, yielding:
\begin{equation}
\sum_{k} (\rho_{i}^{k+1} - \rho_{i}^{k})\,e_{k}(\tilde{t}^{q}) = -\sum_{l}\incidencederivative{1}{0}_{il}\varsigma_{l}^{q}\;.
\label{eq::system_advection_discretization_11}
\end{equation}
Recall that $\rho_{i}^{k}$ is the discrete degree of freedom of the advected quantity at the $t^{k}$ instants of time associated to Gauss-Lobatto nodes and $\varsigma_{l}^{q}$ is the discrete degree of freedom of the fluxes of the advected quantity at the $\tilde{t}^{q}$ instants of time associated to the Gauss nodes, just as stated for the systems of ordinary differential equations.

\subsection{Interior product}

The discretization of the interior product is done using \eqref{eq::adjoint_interior_product}, in the following way:

\begin{svgraybox}
\begin{definition}[\textbf{Discrete interior product}]\label{def::discrete_interior_product}
In one dimension, the discrete interior product $\iota_{\vec{v},h}:\Lambda^{1}_{h}\spacemap\Lambda^{0}_{h}$ is such that:
\begin{equation}
\inner{\iota_{\vec{v},h}\,\kdifform{\alpha}{1}_{h}}{\kdifform{\epsilon}{0}_{i}}_{L^{2}} = \inner{\kdifform{\alpha}{1}_{h}}{\vec{v}^{\flat}\wedge\kdifform{\epsilon}{0}_{i}}_{L^{2}},\quad\forall\kdifform{\epsilon_{i}}{0}\in\Lambda^{0}_{h}\label{eq::discrete_interior_product}
\end{equation}
where $\vec{v}^{\,\flat}=\kdifform{\nu}{1}\in \Lambda^{1}$ and $\kdifform{\alpha}{1}_{h}\in\Lambda^{1}_{h}$.
\end{definition}
\end{svgraybox}
In this way one satisfies the duality pairing between the interior product and the wedge product in the discrete setting. 
%Another point is that an open issue in the advection of differentiable forms lies in the discretization of vector fields. The formulation presented here of the discrete interior product circumvents this by using the $\vec{v}^{\flat}=\nu^{(1)}$ instead of the vector field. This allows one to use the already known discretization of differentiable forms, with its properties. 

%For instance for the advection of $n-forms$ one should use $\vec{v}^{\flat}=\star\nu^{(n-1)}$ since, in this way, one is able to represent correctly the divergence of the advecting vector field, which is a quantity that is responsible for the exact preservation of the advected quantity.

Partitioning the domain $\Omega$ in a spectral element cell complex one can apply the discretization of the interior product in each spectral element, obtaining:
\begin{equation}
\sum_{i}\rho_{i}(t)\inner{\kdifform{\epsilon_{i}}{1}}{\nu^{(1)}\wedge\kdifform{\epsilon_{j}}{0}}_{L^{2}} = \sum_{i}\varsigma_{i}(t)\inner{\kdifform{\epsilon_{i}}{0}}{\kdifform{\epsilon_{j}}{0}}_{L^{2}}\;,\quad\forall\kdifform{\epsilon_{j}}{0}\in \Lambda^{0}_{h}\;.
\label{eq::discrete_lie_derivative}
\end{equation}

\subsection{Putting things together: advection}
The complete discrete systems becomes:
%Now it is possible to setup the algebraic system which represents the advection equation, \eqref{eq::system_advection}. For that one only needs to solve simultaneously the system of equations that represent the evolution in time of the advected quantity, \eqref{eq::system_advection_discretization_11}, and the equations that define the instantaneous fluxes (Lie derivative), \eqref{eq::discrete_lie_derivative}:
\begin{equation}
\begin{cases}
\sum_{k} (\rho_{i}^{k+1} - \rho_{i}^{k})\tilde{e}_{k}(\tilde{t}^{q}) = -\sum_{l}\incidencederivative{1}{0}_{il}\varsigma_{l}^{q}\\
\\
\sum_{i,k}\rho_{i}^{k}\kdifform{\epsilon_{k}}{0}(\tilde{t}^{q})\inner{\kdifform{\epsilon_{i}}{1}}{\star\nu^{(0)}\wedge\kdifform{\epsilon_{j}}{1}}_{L^{2}\Lambda^{1}(\Omega_{m})} = \sum_{i}\varsigma_{i}^{q}\inner{\kdifform{\epsilon_{i}}{0}}{\kdifform{\epsilon_{j}}{0}}_{L^{2}\Lambda^{0}(\Omega_{m})},\,\forall\kdifform{\epsilon_{j}}{0}\in\Lambda^{0}(\Omega_{m})
\end{cases}
\label{eq::system_advection_discretization_02}
\end{equation}

%It is important to note that this results in a staggered scheme in time, where the advected quantity is discretized in Gauss-Lobatto time instants but the flux  is discretized in the Gauss time instants. For this reason the terms $\sum_{k}\rho_{i}^{k}\kdifform{\epsilon_{k}}{0}(\tilde{t}^{q})$, which correspond to interpolation in time, appear. 

%Notice that the equation that discretizes the Lie derivative and therefore the fluxes of the advected quantity is done in the Gauss time instants, $\tilde{t^{q}}$. This means that the discrete fluxes obtained from the discrete Lie derivative are the same as the ones that appear in the time evolution equation. If one had discretized the Lie derivative equation at the Gauss-Lobatto instants, $t^{k}$, the fluxes computed here would have been interpolated from the ones used in the time evolution equation. This would generate an incorrect value of the divergence of the fluxes which would lead to a loss or gain of advected quantity in time. 
 
 \section{Numerical results}
This approach was applied to the two dimensional solution of an advected sine wave and a sine bell in a constant velocity field $\vec{v}=\vec{e}_{x}$: $\rho^{(2)}(x,y) = \sin(\pi x)\sin(\pi y)\ederiv x\ederiv y$ (sine wave) and $\rho^{(2)}(x,y) = \sin(2\pi x)\sin(2\pi y)\ederiv x\ederiv y$ if $(x,y) \in[0,0.5]\times[0,0.5]$ and $\rho^{(2)}(x,y) = 0$ in $(x,y) \in \mathbb{R}^{2}\backslash[0,0.5]\times[0,0.5]$  (sine bell), on a domain with periodic boundary conditions.

\begin{figure}[ht]
\centering
\subfigure{
\includegraphics[width=0.4\textwidth]{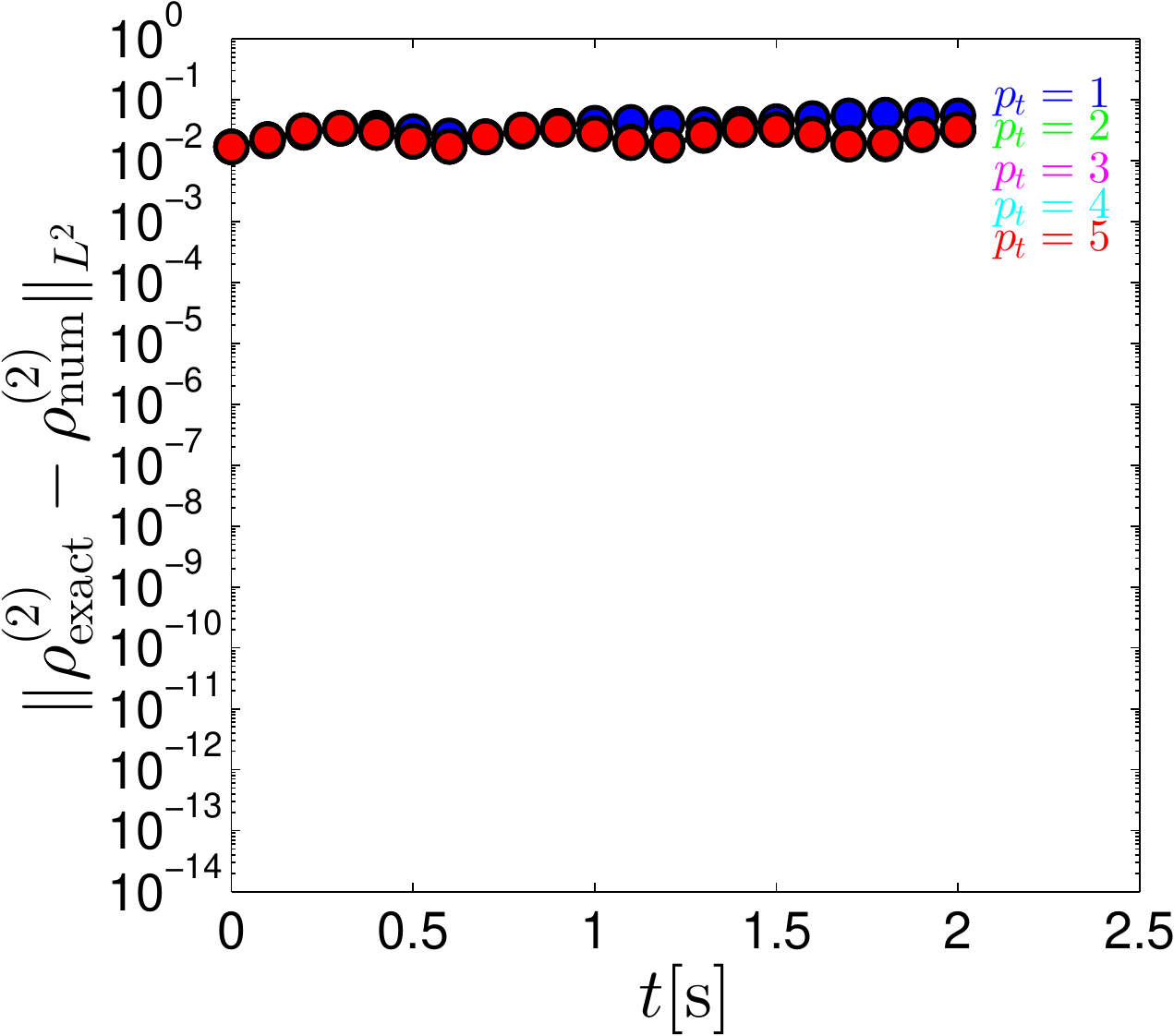}
\label{fig:advection_time_error_p3}
}
\subfigure{
\includegraphics[width=0.4\textwidth]{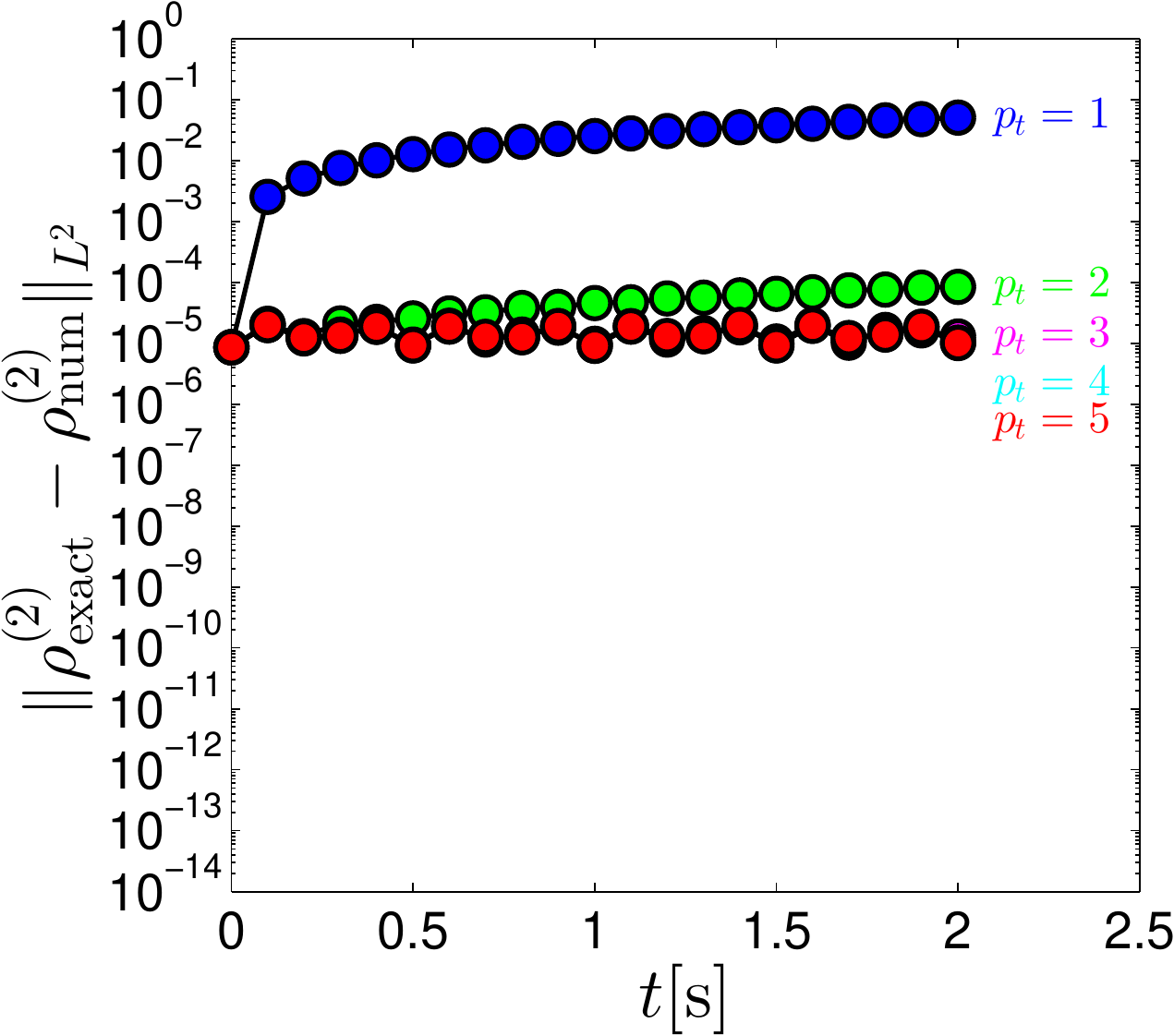}
\label{fig:advection_time_error_p6}
}
\subfigure{
\includegraphics[width=0.4\textwidth]{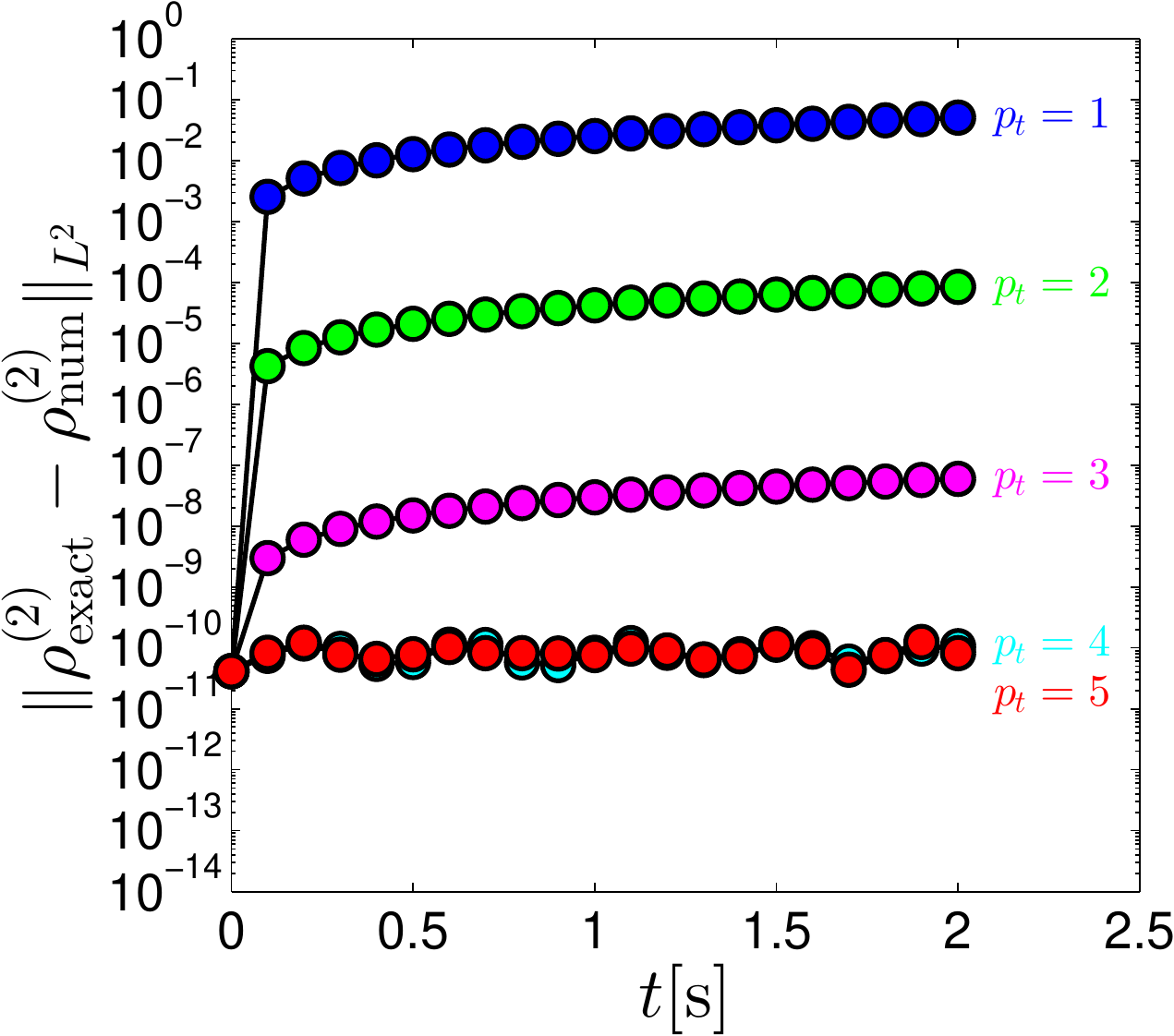}
\label{fig:advection_time_error_p10}
}
\subfigure{
\includegraphics[width=0.4\textwidth]{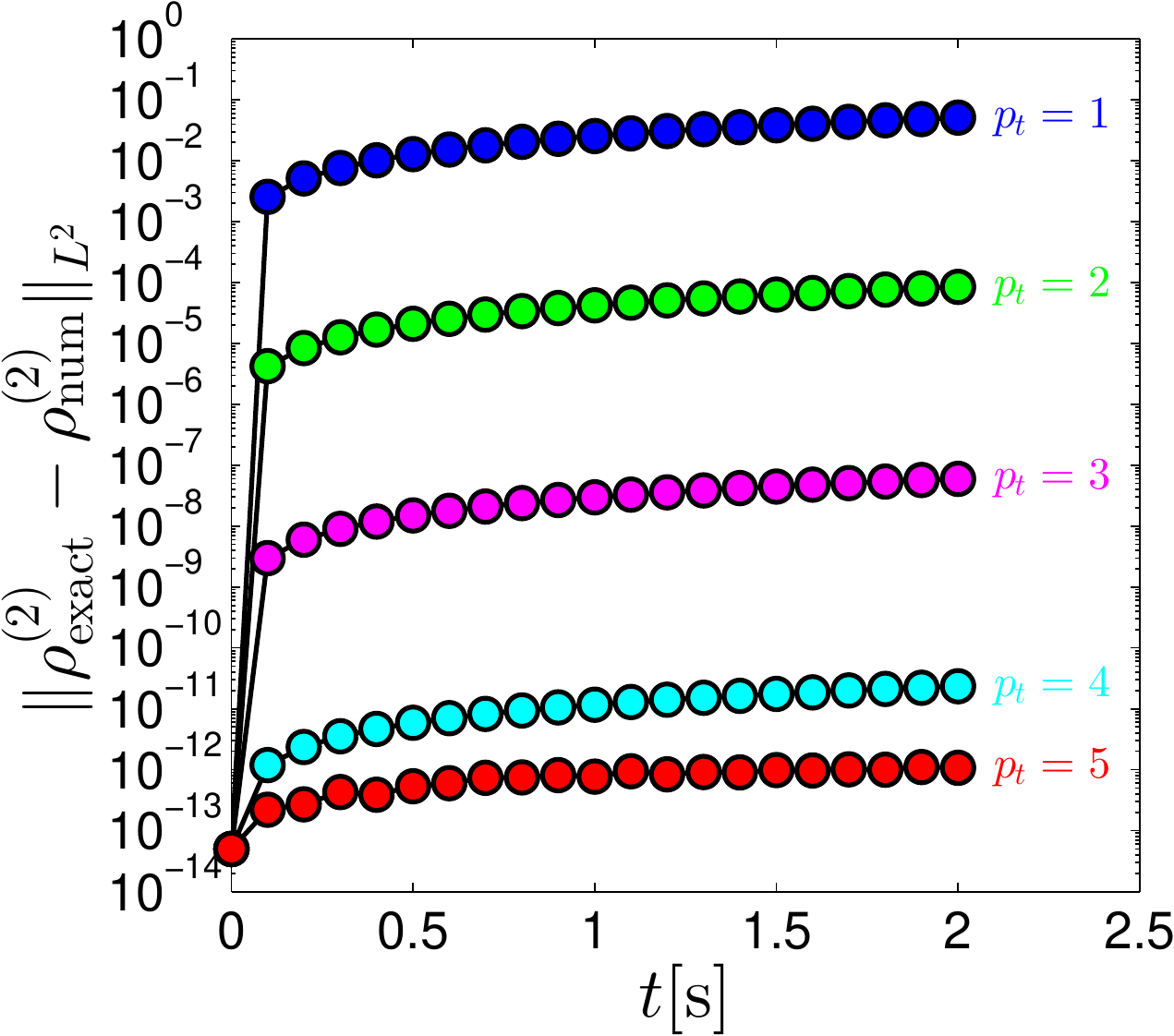}
\label{fig:advection_time_error_p12}
}
\caption{Error in time of the numerical solution of \eqref{eq::system_advection}  with $\vec{v}=\vec{e}_{x}$ and $4\times4$ elements of order $p=3$, $p=6$, $p=10$ and $p=12$ (from left to right and top to bottom) and $\Delta t = 0.1s$, for the sine wave $\rho^{(2)}(x,y) = \sin(\pi x)\sin(\pi y)\ederiv x\ederiv y$. As shown, the error in the solution increases with time due to the inaccuracy of time integration. When time integration is accurate enough the error in the initial state is preserved. Here $p_{t}$ denotes the polynomial degree in time.}
\label{fig:advection_time_error_p3-12}
\end{figure}

In \figref{fig:advection_time_error_p3-12} the error in time of the numerical solution of \eqref{eq::system_advection} for a mesh of $4\times4$ elements with a $\Delta t = 0.1s $ and various polynomial orders in space, $p$, and time, $p_{t}$, is presented. The initial error, due to the discretization, is conserved, as long as the time integration is sufficiently accurate.

\begin{figure}[ht]
\centering
\subfigure{
\includegraphics[width=0.4\textwidth]{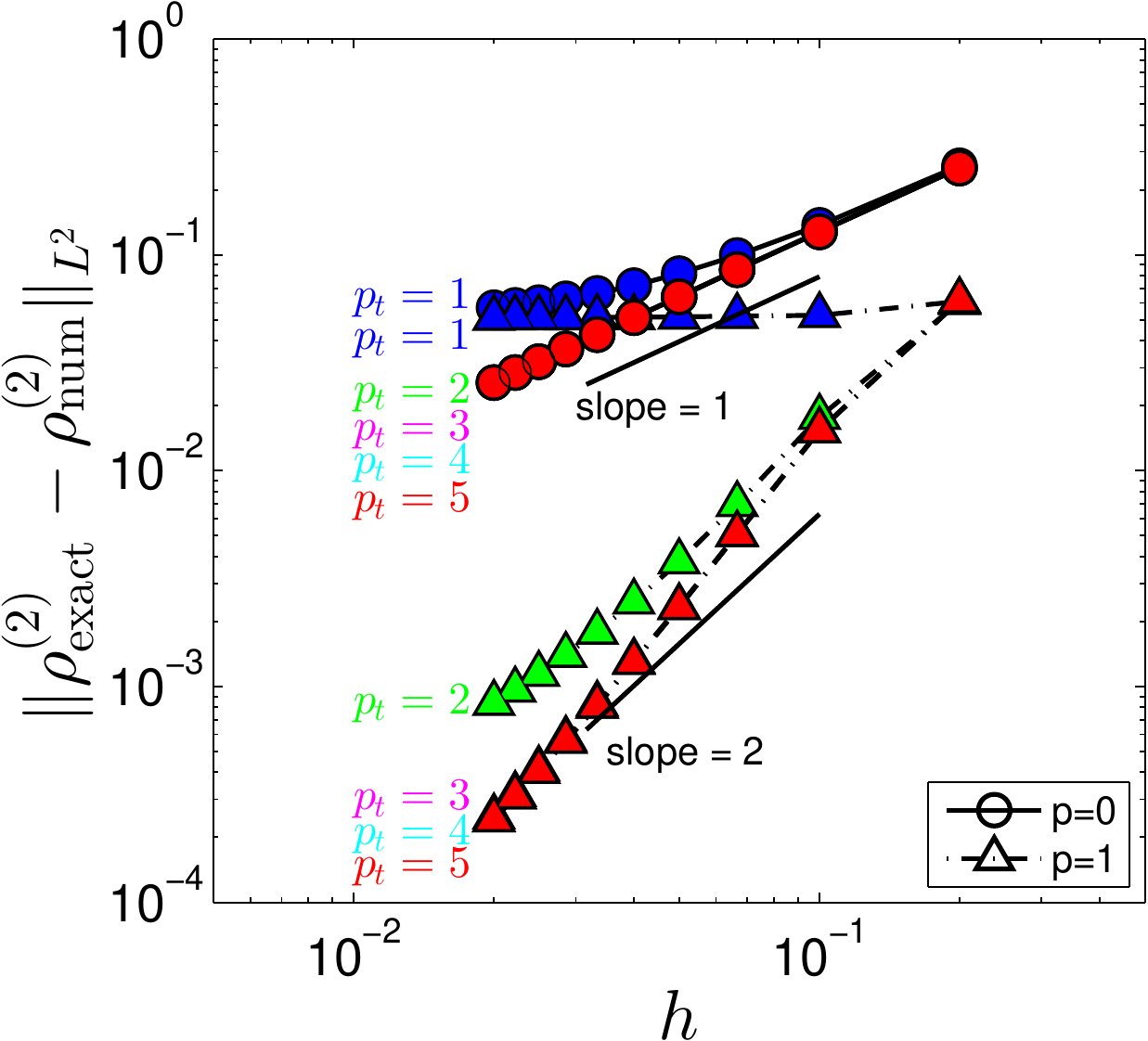}
\label{fig:advection_h_convergence}
}
\subfigure{
\includegraphics[width=0.4\textwidth]{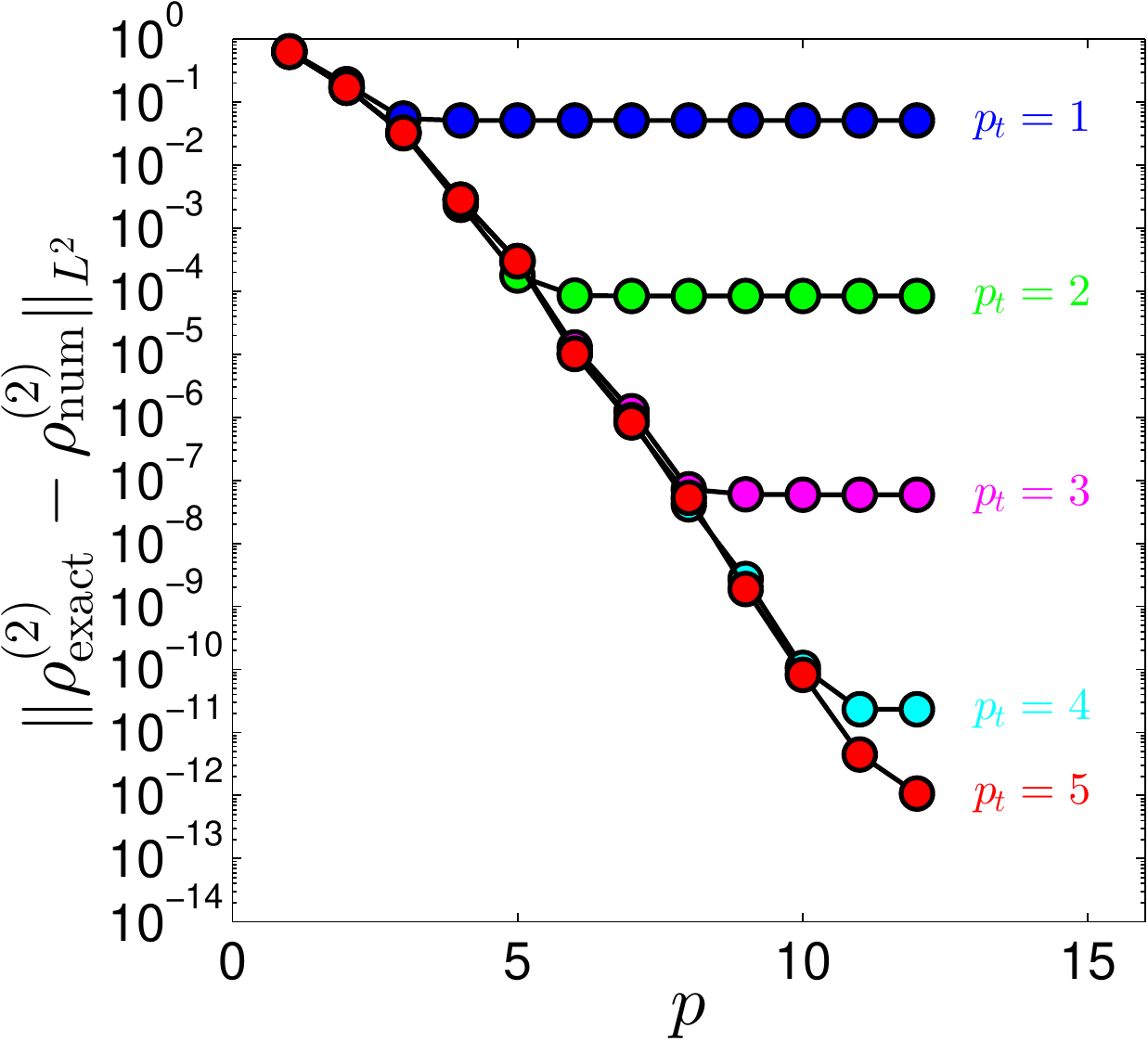}
\label{fig:advection_p_convergence}
}
\caption{Left: $h$ convergence in space for the advection of a sine wave with $\Delta t =0.1s$. Right: $p$ convergence in space for the advection of a sine wave, $4\times 4$ elements and $\Delta t =0.1s$.}
\label{fig:advection_space_error_convergence}
\end{figure}

In \figref{fig:advection_space_error_convergence}, the $h$- and $p$-convergence plots are shown for different values of the order of the time integration scheme, $p_{t}$, and $\Delta t=0.1s$. It is possible to see that the method presents algebraic $h$-convergence rates of order $(p+1)$ as long as the time integration error does not dominate the spatial one. The method shows a spectral $p$-convergence as soon as the time integration is accurate enough.

\begin{figure}[ht]
\centering
\includegraphics[width=0.4\textwidth]{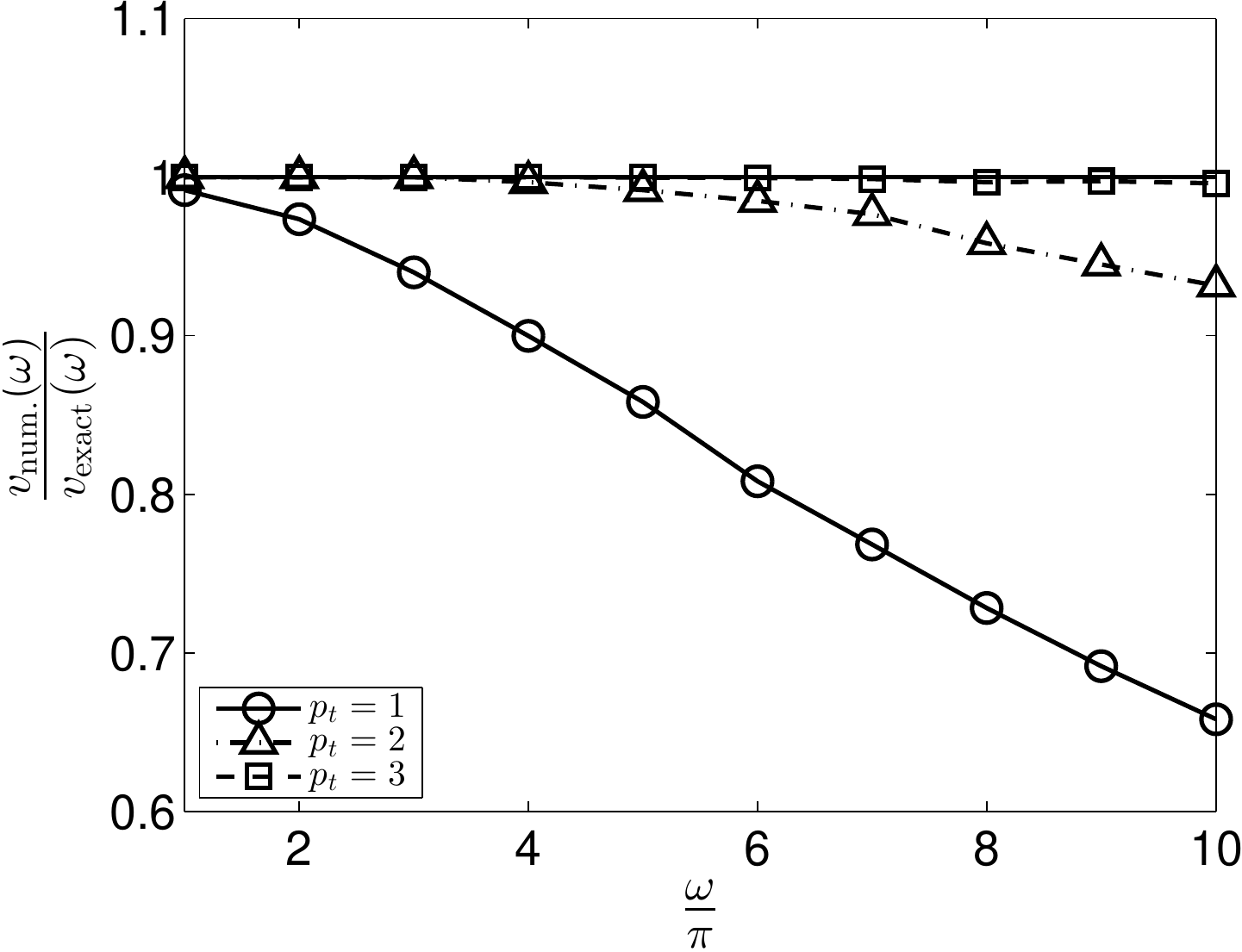}
\caption{Error in velocity as a function of the frequency of the advected sine wave: numerical dispersion. $p=10$, $\Delta t=0.1s$ and $n=4\times 4$ elements.}
\label{fig:advection_numerical_dispersion}
\end{figure}

In \figref{fig:advection_numerical_dispersion} the error on the velocity is presented as a function of the advected sine wave frequency. This figure shows that the numerical method introduces an artificial dispersion if the time scheme is not accurate enough. 

\begin{figure}[ht]
\centering
\includegraphics[width=0.4\textwidth]{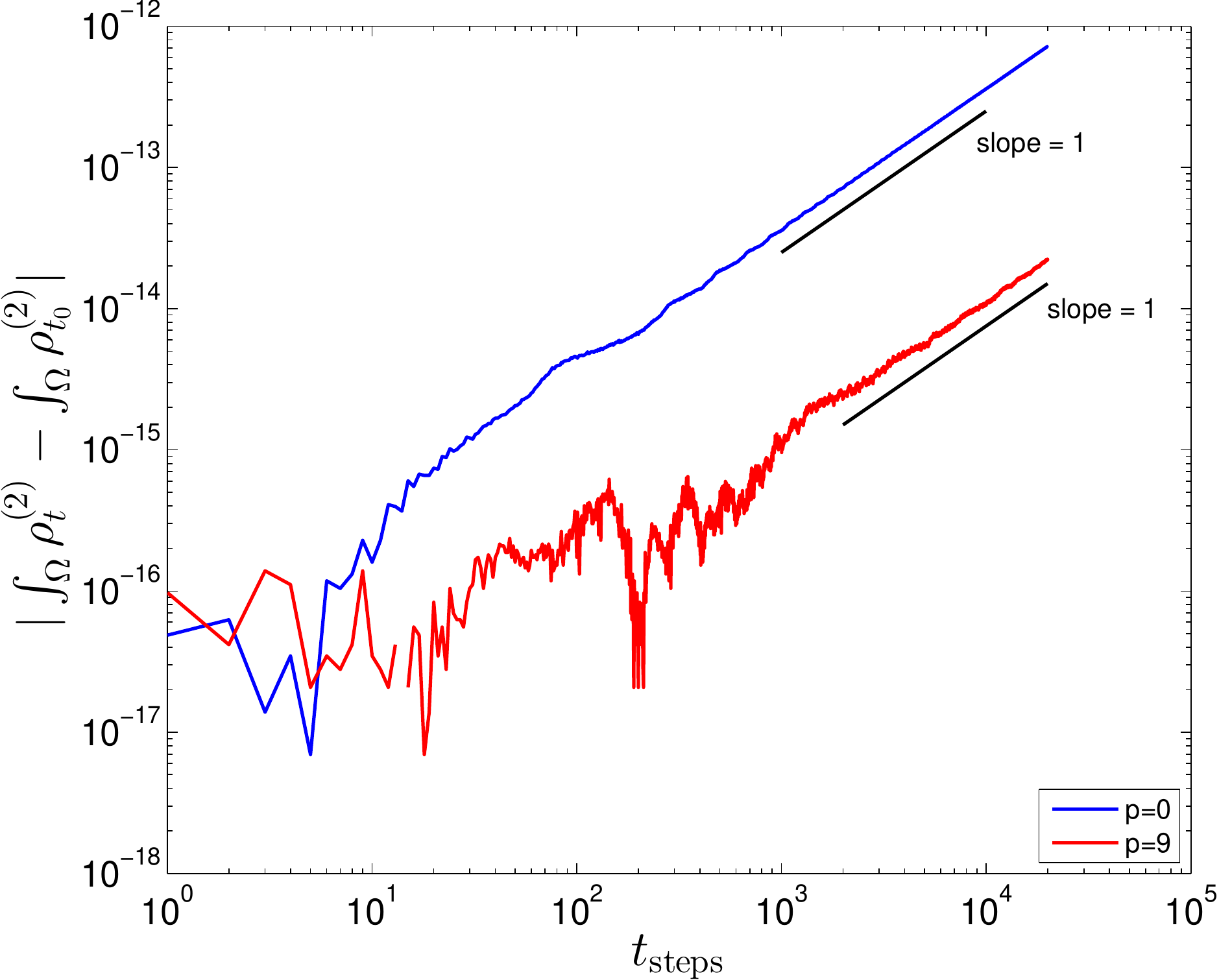}
\caption{Sum of the local errors of the advected 2-form for a sine bell in a velocity field $\vec{v}=\vec{e}_{x}$, with $50\times 50$ elements of order $p=0$ (blue) and $4\times 4$ elements of order $p=9$ (red), $\Delta t = 0.01$ and $p_{t}=2$.}
\label{fig:advection_mass_error_sine_bell}
\end{figure}

Another fundamental aspect is the conservation of the advected quantity. \figref{fig:advection_mass_error_sine_bell} shows the mass error in time, that is: $\int_{\Omega}\kdifform{\rho}{2}_{t}-\int_{\Omega}\kdifform{\rho}{2}_{t_{0}}$. The error goes from the zero machine in the first $10^{3}$ time steps while thereafter it steadily increases. Notice that even after $2\times10^{4}$ time steps the error is still below $10^{-12}$.

% In \figref{fig:advection_mass_square_error_sine_bell} is shown the square of mass error in time, that is: $\int_{\Omega}\kdifform{\rho}{2}_{t}\wedge\star\kdifform{\rho}{2}_{t}-\int_{\Omega}\kdifform{\rho}{2}_{t_{0}}\wedge\star\kdifform{\rho}{2}_{t_{0}}$. As can be seen the error for the case $p=0$ is of the order of machine accuracy, hence the advected quantity is conserved in time. For the case where $p=10$ we see the error bounded between $10^{-4}$ and $10^{-9}$. A closer look shows that the error is minimal, $10^{-9}$ at intervals with a period of 0.5s. This value of 0.5s corresponds to one element width. Hence the error in the integral of the square of the advected quantity appears when the sine bell is crossing elements and returns to the initial value periodically. This happens due to the fact that when the sine bell is between two (or more) elements bigger discontinuities in the boundary of the elements arise.

\begin{figure}[ht]
\centering
\subfigure{
\includegraphics[width=0.3\textwidth]{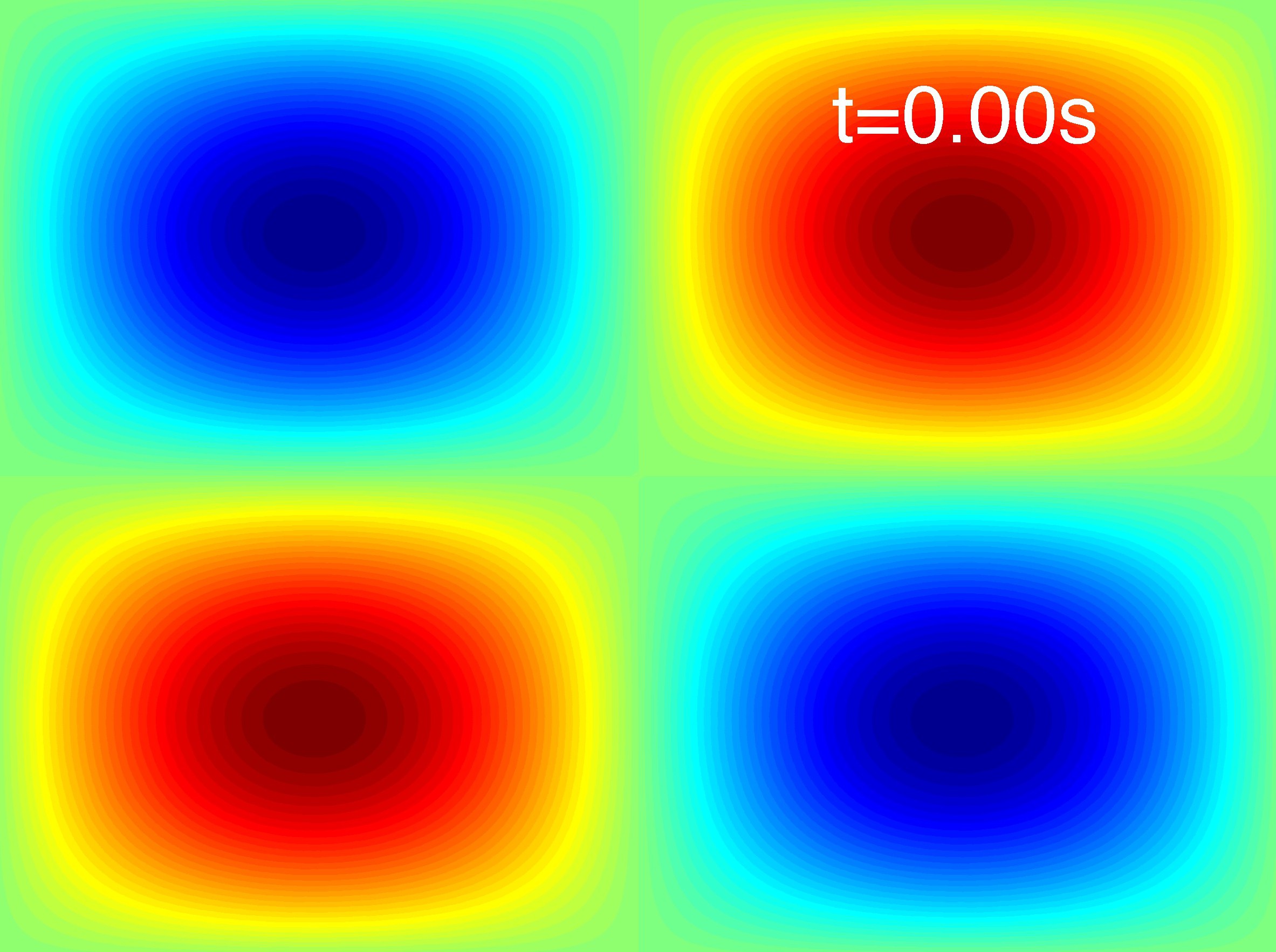}
\label{fig:advection_reversibility_00s}
}
\subfigure{
\includegraphics[width=0.3\textwidth]{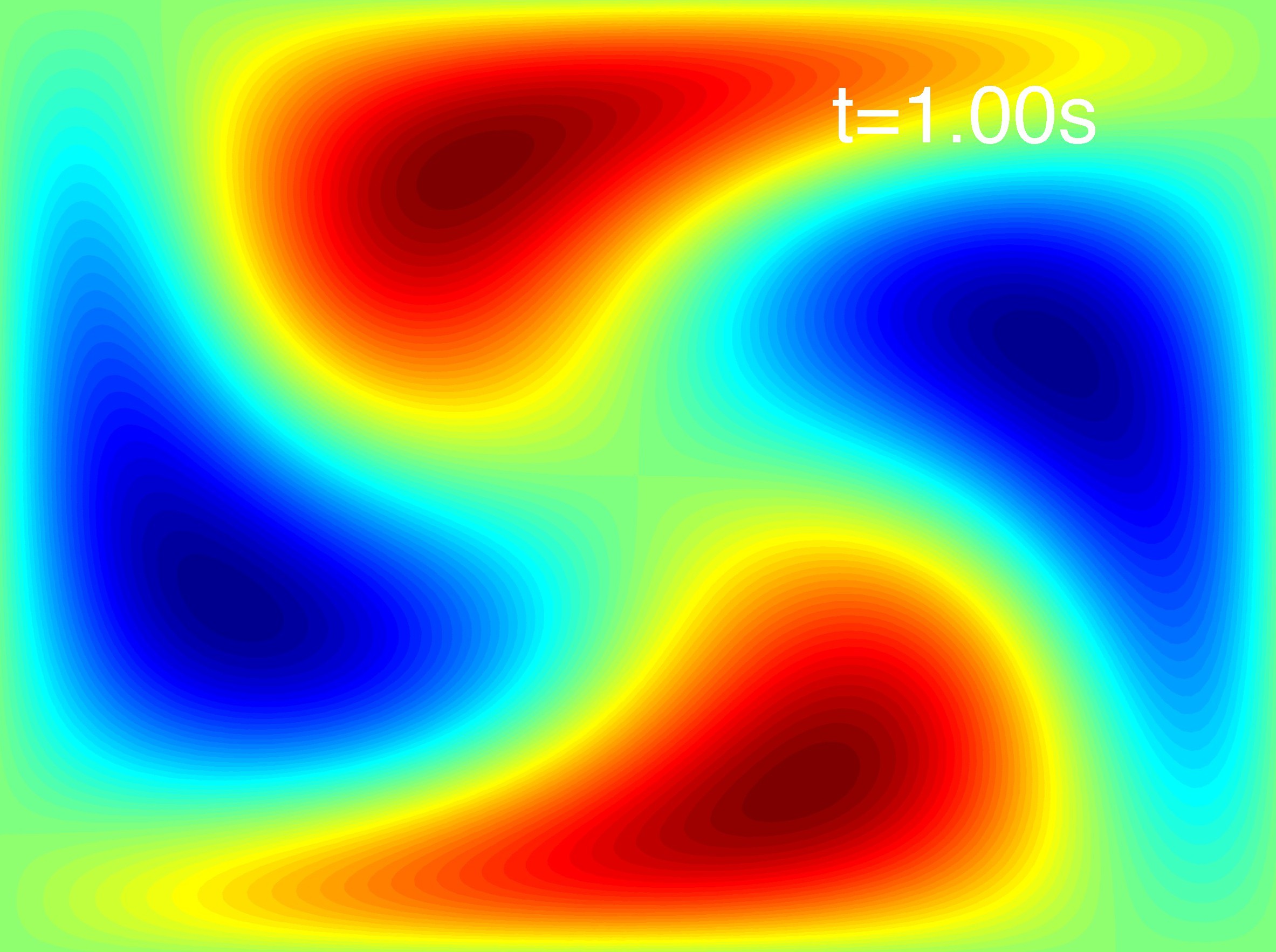}
\label{fig:advection_reversibility_10s}
}
\subfigure{
\includegraphics[width=0.3\textwidth]{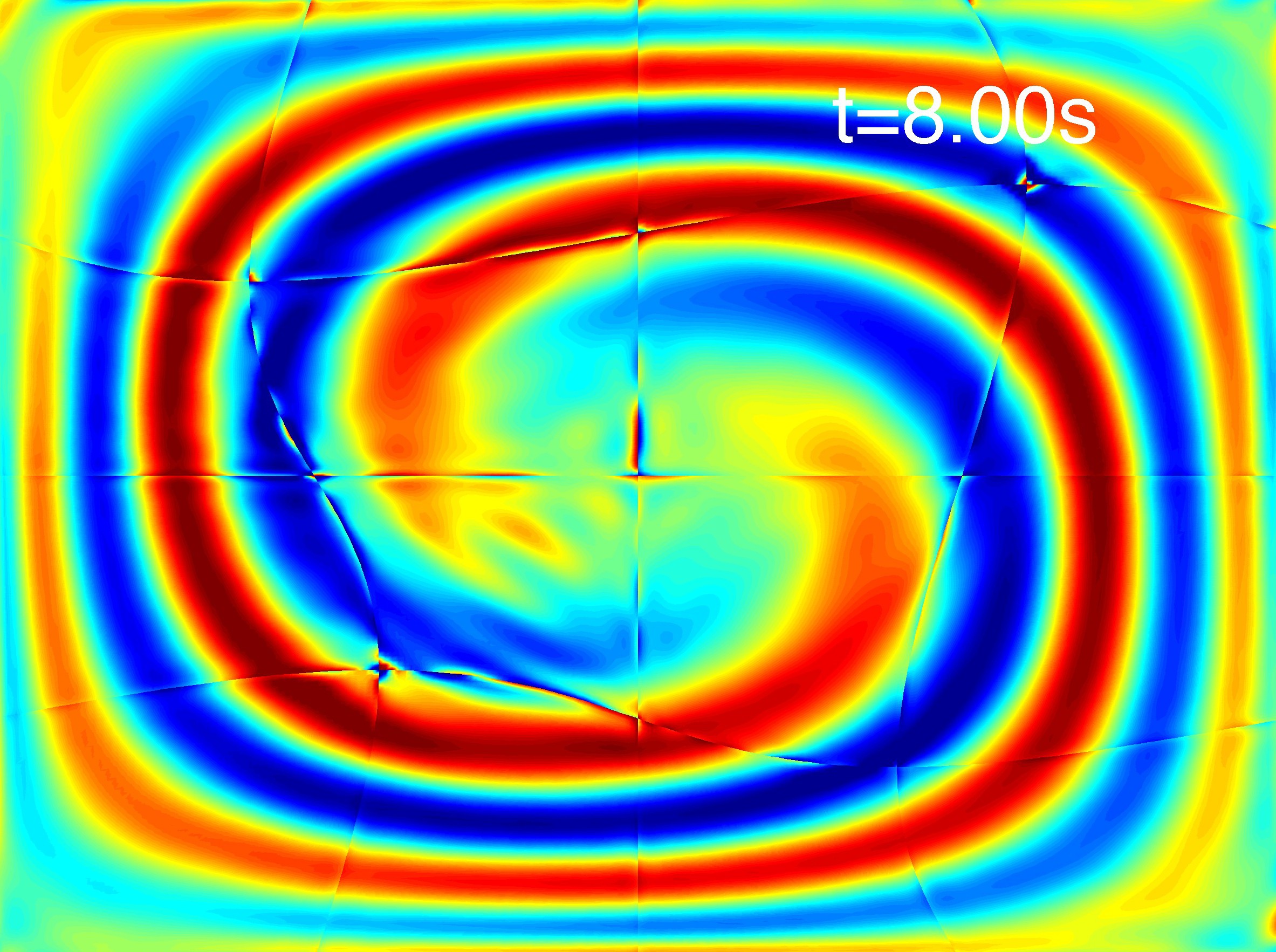}
\label{fig:advection_reversibility_80s}
}
\subfigure{
\includegraphics[width=0.3\textwidth]{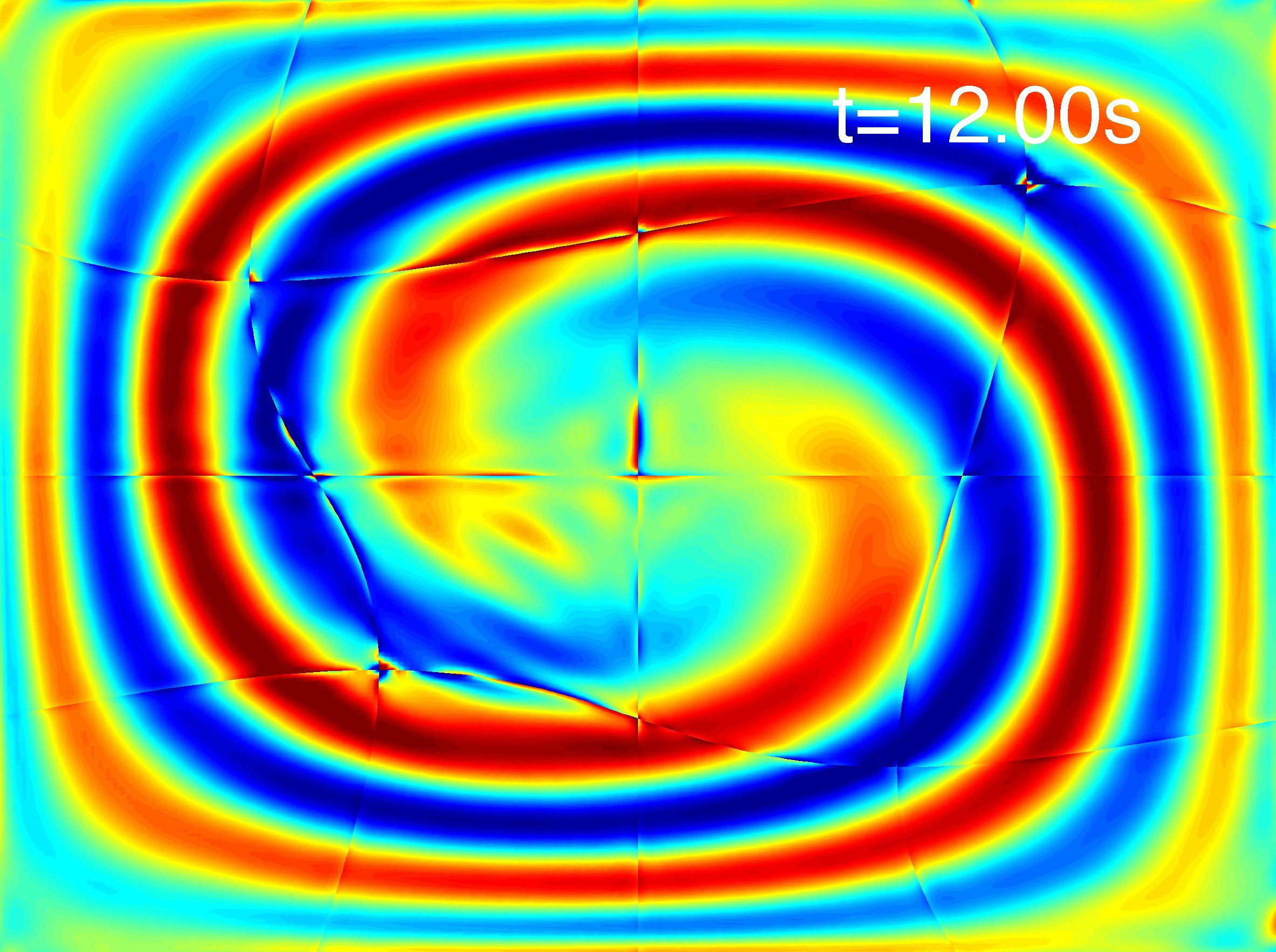}
\label{fig:advection_reversibility_120s}
}
\subfigure{
\includegraphics[width=0.3\textwidth]{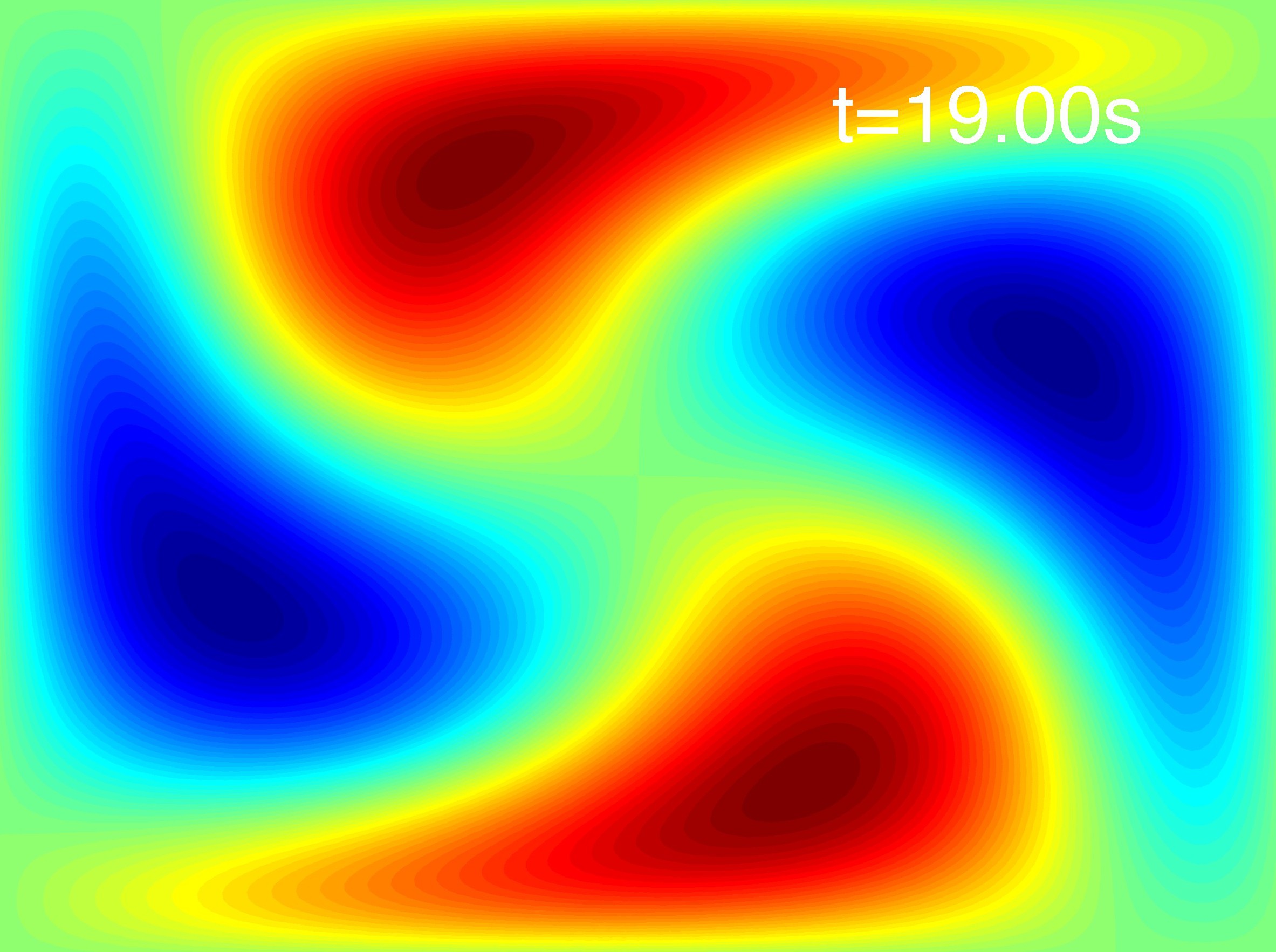}
\label{fig:advection_reversibility_190s}
}
\subfigure{
\includegraphics[width=0.3\textwidth]{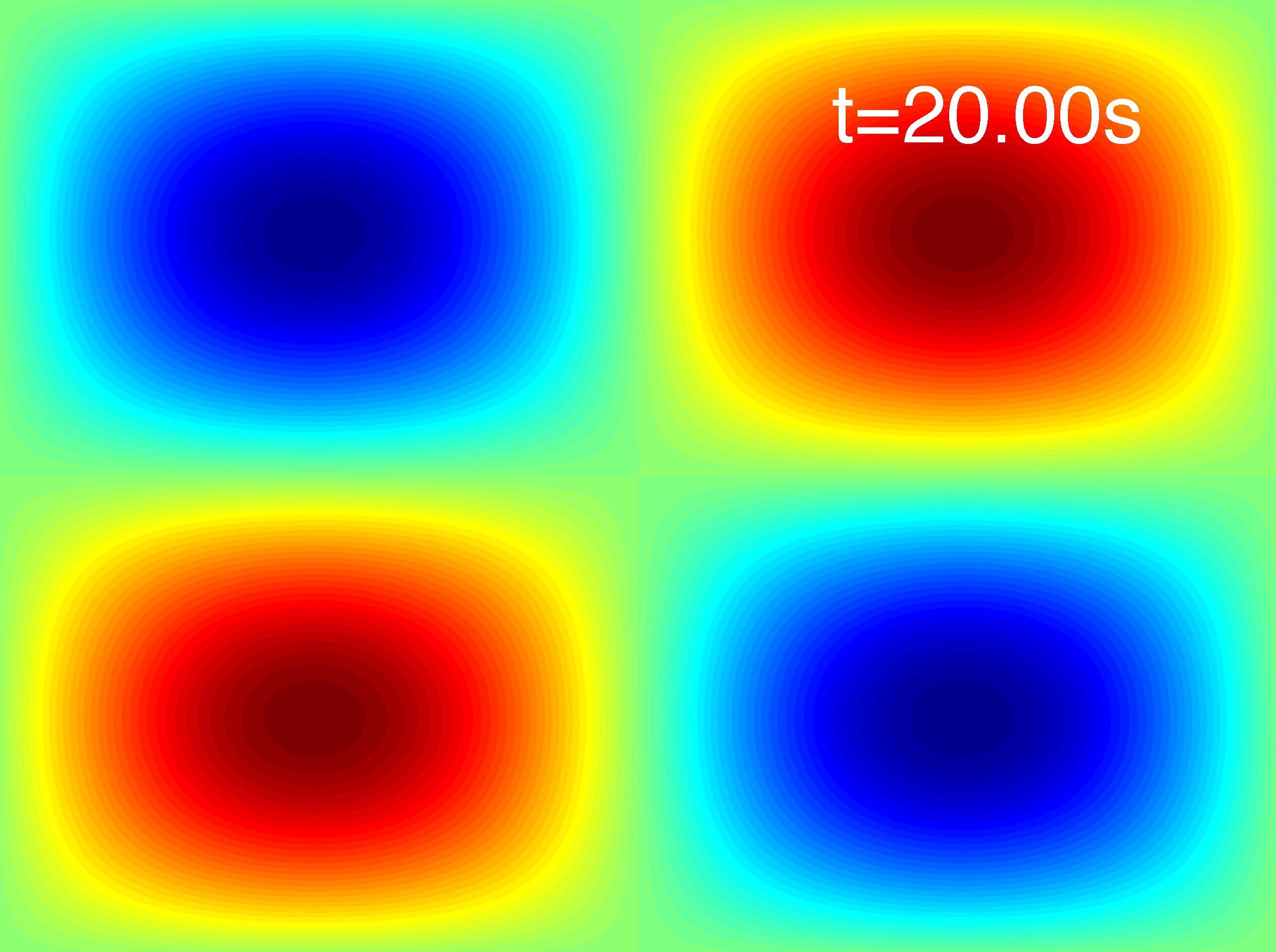}
\label{fig:advection_reversibility_200s}
}
\caption{From left to right and from top to bottom: advection of a sine wave of frequency $\omega=\pi$ on a Rudman vortex with $4\times 4$ curved elements of order $p=9$, $\Delta t=0.1s$ and time integration of order $p_{t}=2$, on a distorted mesh. At time $t=10$, the flow field is reversed. At times $t=8.0s$ and $t=12.0s$ the mesh is visible in the solution. See \texttt{http://www.youtube.com/watch?v=QmoJyqtk9YA} for animation.}
\label{fig:advection_rudman_vortex}
\end{figure}

In \figref{fig:advection_rudman_vortex} one can see the advection of a sine wave of frequency $\omega=\pi$ in a Rudman vortex for 100 time steps after which the direction of the flow is reversed and the calculation is continued for another 100 time steps, with $4\times 4$ curved elements of order $p=9$, $\Delta t=0.1s$ and time integration of order $p_{t}=2$, on a distorted mesh. The mimetic advection enables one to recover the initial solution, thus demonstrating that the integration method is reversible.

% \begin{figure}[ht]
% \centering
% \subfigure{
% \includegraphics[width=0.4\textwidth]{./testCases/advection/results/sine_bell_mass_square_error_time-crop.pdf}
% \label{fig:sine_bell_mass_conservation_p1}
% }
% \subfigure{
% \includegraphics[width=0.4\textwidth]{./testCases/advection/results/sine_bell_mass_square_error_time_p10-crop.pdf}
% \label{fig:sine_bell_mass_conservation_p10}
% }
% \caption{Error in the integral of the square of advected 2-form for the advection of a sine bell on a velocity field $\vec{v}=\vec{e}_{x}$, with $50\times 50$ elements of order $p=0$ (left) and $4\times 4$ elements of order $p=10$ (right), $\Delta t = 0.01$ and $pTime=2$.}
% \label{fig:advection_mass_square_error_sine_bell} 
% \end{figure}

%\section{SUMMARY AND OUTLOOK}
%\label{Section::Summary}
%
%In this paper we presented a first extension of the framework \cite{kreeft2011mimetic} to the numerical solution of the linear advection equation. In terms of L$^{2}$-error it was observed that the rate of convergence $h-$ refinement  in space is $p+$ and for $p-$refinement in space one can observe spectral convergence, as long as not limited by the order of time integration. Furthermore, artificial dispersion is observed, which decreases with the order of time integration $p_{t}$. One major issue in our formulation lies on the discontinuity between elements. Current work focuses on the solution of Euler's equation for incompressible flow by applying this approach to both the non-linear advection of vorticity and of momentum.

\begin{acknowledgement}
The authors would like to thank the valuable comments of both reviewers and the funding received by FCT - Foundation for science and technology Portugal through SRF/BD/36093/2007 and SFRH/BD/79866/2011.
\end{acknowledgement}

\bibliographystyle{plain}
\bibliography{./AdvectionICOSHAM2012Submit_reviewed.bib}

\end{document}